\documentclass{amsart}
\usepackage{amsmath,amssymb,latexsym,mathrsfs}
\usepackage[mathscr]{eucal}
\usepackage{enumerate}
\usepackage{color}

\theoremstyle{remark}

\begin{document}
\title[Monotonic Extensions of Lipschitz Maps]{Order-Preserving Extensions
of Hadamard Space-Valued Lipschitz Maps}
\author{Edoardo Gargiulo}
\address{Courant Institute of Mathematical Sciences, New York University}
\email{eg3504@nyu.edu}
\author{Efe A. Ok}
\address{Department of Economics and Courant Institute of Mathematical
Sciences, New York University}
\email{efe.ok@nyu.edu}
\date{January 25, 2026}
\thanks{This paper is in final form and no version of it will be submitted
for publication elsewhere.}

\begin{abstract}
We study the problem of extending any order-preserving Lipschitz function
that maps a subset of a partially ordered Hilbert space $X$ into a Hadamard
poset $Y$ without increasing its Lipschitz constant and preserving its
monotonicity. This sort of an extension is always possible when $X$ is
one-dimensional. However, when $\dim X\geq 2$ and $Y$ satisfies some fairly
weak conditions, it holds (universally) if and only if the order of $X$ is
trivial. The conditions on $Y$ are satisfied by any Hilbert poset.
Therefore, as a special case of our main result, we find that there is no
order-theoretic generalization of Kirszbraun's theorem.
\end{abstract}

\maketitle

\section{Introduction}

A pair $(X,Y)$ of metric spaces has the \textit{Lipschitz extension property}
if for any nonempty $S\subseteq X$ and any 1-Lipschitz map $f:S\rightarrow
Y, $ there exists a 1-Lipschitz $F:X\rightarrow Y$ that extends $f.$ One of
the central problems in Lipschitz analysis is to determine those $(X,Y)$
that has the Lipschitz extension property.\footnote{%
A stellar reference on this problem and its variants is Benyamini and
Lindenstrauss \cite{B-L}$.$} Two classical results delineate the scope of
this problem. The McShane-Whitney theorem (of \cite{McS} and \cite{Whit})
asserts that $(X,\mathbb{R})$ has the Lipschitz extension property for every
metric space $X,$ while Kirszbraun's theorem (of \cite{Kirszbraun}) says $%
(X,Y)$ has the Lipschitz extension property for every Hilbert spaces $X$ and 
$Y$.\footnote{%
To be precise, Kirszbraun proved this for maps between two
finite-dimensional Euclidean spaces. Valentine \cite{Valentine} noticed that
the theorem extends to the context of Hilbert spaces without modification.
Some authors thus refer to this result as the Kirszbraun-Valentine theorem.}

When metric spaces carry additional structure, it is natural to ask whether
Lipschitz extensions can be obtained so as to preserve that structure. In
particular, when the domain and the codomain are partially ordered metric
spaces, one may require extensions to be order-preserving in addition to
being Lipschitz. To be precise, we say that a pair $(X,Y)$ of partially
ordered metric spaces has the \textit{monotone Lipschitz extension property}
if for any nonempty $S\subseteq X$ and order-preserving 1-Lipschitz map $%
f:S\rightarrow Y,$ there is an order-preserving 1-Lipschitz $F:X\rightarrow
Y $ with $F|_{S}=f.\footnote{%
The continuous version of this problem for real-valued maps (i.e., extending
a monotonic continuous real function on a subset of a topological poset to
the entire space monotonically and continuously) has received abundant
attention. The literature on this problem has started with Nachbin \cite{Na}%
, and found profound applications in decision theory, mathematical
economics, and optimal transport theory. For recent contributions to this
literature, see \cite{Yamazaki}, and references cited therein.}$

Ok \cite{Ok} has recently characterized the class of all partial orders on $%
X $ relative to which $(X,\mathbb{R})$ has the monotone Lipschitz extension
property. As this class contains the equality relation on $X$ (which is
unable to compare any two distinct points in the space), that result
provides an order-theoretic generalization of the McShane-Whitney theorem.

The purpose of the present paper is to investigate to what extent the
coverage of Kirszbraun's theorem can be extended in this manner. We are thus
particularly interested in the case where $X$ and $Y$ are Hilbert posets (=
Hilbert spaces vector ordered by means of pointed closed convex cones), but
our main result not really linear, and it applies to a much richer class of
metric spaces. We begin by observing that a positive result is readily
obtained by affine interpolation when $X=\mathbb{R}$ and $Y$ is any Banach
poset (Theorem 2). When the dimension of $X$ exceeds 1, however, things
change radically.

In the context of Hilbert spaces, the main result of this paper reads as:

{\small \medskip }

{\small \medskip }

\noindent \textsc{The No-Extension Theorem.} \textsl{Let }$X=(X,\succcurlyeq
_{X})$\textsl{\ be a partially ordered Hilbert space with }$\dim X\geq 2$%
\textsl{, and }$Y=(Y,\succcurlyeq _{Y})$\textsl{\ any Hilbert poset. Then, }$%
(X,Y)$\textsl{\ has the monotone Lipschitz extension property if, and only
if, }$\succcurlyeq _{X}$\textsl{\ is trivial.}

{\small \medskip }

{\small \medskip }

\noindent When $\succcurlyeq _{X}$ is trivial (that is, it is the equality
relation), every map on $X$ is automatically order-preserving. Thus, in that
case, the monotone Lipschitz extension property is identical to the
Lipschitz extension property, so all reduces to Kirszbraun's theorem. The
no-extension theorem thus shows that there is no order-theoretic
generalization of Kirszbraun's theorem.

We obtain the no-extension theorem as a special case of a more general
no-extension theorem in which $Y$ belongs to a fairly large class of
Hadamard posets. (A Hadamard poset is a complete CAT(0)-space endowed with a
geodesically hereditary closed partial order; see Section 2.4.) At a
heuristic level, the proof of this latter result rests on a simple, but
rigid, obstruction mechanism that combines order structure of the domain
with asymptotic geometry of the codomain. The key observation is that if
order-preserving 1-Lipschitz extensions into a given Hadamard poset $Y$ were
universally possible, then the partial order on the domain $X$ would
necessarily satisfy a strong metric compatibility condition, called \textit{%
radiality} (Section 2.2). This necessity is detected by a test map on a
two-point set whose order-preserving 1-Lipschitz extension forces a metric
inequality via composition with a monotone Busemann function associated with
an order-preserving geodesic ray in $Y$ (Theorem 3). Thus, to obtain this
necessity, we presume that $Y$ possesses at least one monotone Busemann
function associated with an order-preserving geodesic ray. This condition is
met by any Hilbert poset (Proposition 10), among many other Hadamard posets
(Examples 1-3).

It remains to explore how demanding radiality is. This property is fairly
forgiving in the context of finite metric posets (cf. \cite{Ok}). By
contrast, we find here that it is extremely restrictive on connected metric
posets (Theorem 4). In fact, if a metric space contains a copy of a close
simple loop, such as a normed space with dimension $\geq 2,$ then it cannot
be nontrivially radially ordered (Theorem 5). Putting these observations
together yields a general result (Theorem 8), a special case of which is the
necessity part of the no-extension theorem we stated above.

We conclude the paper by outlining a research agenda for studying the
monotone Lipschitz extension problem quantitatively.

\section{Preliminaries}

\subsection{Posets}

Let $X$ be a nonempty set. A \textit{preorder} on $X$ is a reflexive and
transitive binary relation on $X,$ while a \textit{partial order} on $X$ is
an antisymmetric preorder on $X.$ We refer to the ordered pair $%
(X,\succcurlyeq )$ as a \textit{poset }if $\succcurlyeq $ is a partial order
on $X.$ Obviously, the smallest partial order on $X$ is $\{(x,x):x\in X\};$
we refer to this partial order as the \textit{trivial order} on $X$. A
preorder $\succcurlyeq $ on $X$ is \textit{total} if any two elements $x$
and $y$ of $X$ are $\succcurlyeq $\textit{-comparable}, that is, either $%
x\succcurlyeq y$ or $y\succcurlyeq x$ holds. A total partial order on $X$ is
said to be a \textit{linear order} on $X$.

In what follows, unless otherwise is explicitly stated, we regard $\mathbb{R}%
^{n}$ as a poset relative to the usual coordinatewise order. Obviously, this
order is total iff $n=1$.

Given any poset $(X,\succcurlyeq ),$ we denote the asymmetric part of $%
\succcurlyeq $ by $\succ $, that is, $x\succ y$ means $y\neq x\succcurlyeq
y. $ We also define the binary relation $\succcurlyeq ^{\bullet }$ on $X$ as%
\begin{equation*}
x\succcurlyeq ^{\bullet }y\hspace{0.2in}\text{iff\hspace{0.2in}not }%
y\succcurlyeq x.
\end{equation*}%
Thus $x\succcurlyeq ^{\bullet }y$ means that either $x\succ y$, or $x$ and $%
y $ are not $\succcurlyeq $-comparable. It is plain that $\succcurlyeq
^{\bullet }$ is an irreflexive relation. In general, this relation is
neither symmetric nor asymmetric, nor it is transitive. When $\succcurlyeq $
is total, however, $\succcurlyeq ^{\bullet }$ equals $\succ $.

A function $f:X\rightarrow Y$ from a poset $X=(X,\succcurlyeq _{X})$ to a
poset $Y=(Y,\succcurlyeq _{Y})$ is said to be \textit{order-preserving} (or 
\textit{isotonic}) if for every $x,y\in X$, 
\begin{equation*}
x\succcurlyeq _{X}y\hspace{0.2in}\text{implies\hspace{0.2in}}%
f(x)\succcurlyeq _{Y}f(y).
\end{equation*}%
If $Y=(\mathbb{R},\geq ),$ where $\geq $ is the usual order, we refer to $f$
simply as $\succcurlyeq $\textit{-increasing.}

\subsection{Partially Ordered Metric Spaces}

Let $X=(X,d)$ be a metric space. We say that $(X,\succcurlyeq )$ is a 
\textit{partially ordered metric space} if $\succcurlyeq $ is a partial
order on $X.$ If, in addition, $\succcurlyeq $ is a closed subset of $%
X\times X,$ we refer to $(X,\succcurlyeq )$ as a \textit{metric poset.}

By a \textit{radially ordered metric space}, we mean a partially ordered
metric space $(X,\succcurlyeq )$ such that 
\begin{equation}
x\succcurlyeq ^{\bullet }y\succ z\hspace{0.2in}\text{implies\hspace{0.2in}}%
d(x,z)\geq d(x,y)  \label{rd1}
\end{equation}%
and%
\begin{equation}
x\succ y\succcurlyeq ^{\bullet }z\hspace{0.2in}\text{implies\hspace{0.2in}}%
d(x,z)\geq d(y,z)\text{.}  \label{rd2}
\end{equation}%
(In this case, we refer to the partial order $\succcurlyeq $ on $X$ as 
\textit{radial}.) This concept will play an essential role in our analysis
of the monotonic Lipschitz extension problem.

The trivial order on $X$ is obviously radial. \cite{Ok} provides several
other examples of radial orders on a metric space.

{\small \medskip }

\noindent \textsc{Remark.} Radial partial orders were introduced in Ok \cite%
{Ok} who showed that the property of radiality is necessary and sufficient
for all order-preserving 1-Lipschitz real-valued maps on a subset of a
partially ordered metric space to be extendable to such maps on the entire
space. As such, it is not surprising that it will play an important role for
the present work, but the role of radiality here will be a bit more nuanced.
It may be worth noting that every radially ordered metric space $%
(X,\succcurlyeq )$ is \textit{radially convex}, that is, $d(x,z)\geq \max
\{d(x,y),d(y,z)\}$ holds whenever $x\succ y\succ z.\footnote{%
Radially convex metric posets are studied in topological order theory (cf. 
\cite{B-O}, \cite{Carruth} and \cite{Ward}, among many others). They are
also used in the topological analysis of smooth dendroids; (cf. \cite{Fugate}%
).}$ The converse is false, but a linearly ordered metric space is radial if
and only if it is radially convex.

{\small \medskip }

It may come as a small surprise that every radial partial order is closed.
This was proved in \cite{Ok} in a very indirect manner. As we will need it
below, we give an elementary proof of this fact.

{\small \medskip }

{\small \medskip }

\noindent \textsc{Lemma 1.} \textsl{Every radially ordered metric space }$%
(X,\succcurlyeq )$ \textsl{is a metric poset.}

{\small \medskip }

\noindent \textit{Proof.} Take any $x,x_{1},x_{2},...,y,y_{1},y_{2},...\in X$
such that $x_{m}\succcurlyeq y_{m},$ $x_{m}\rightarrow x,$ and $%
y_{m}\rightarrow y.$ To derive a contradiction, suppose $x\succcurlyeq y$ is
false. Let us distinguish between two cases. First, suppose $x\succcurlyeq
y_{m}$ for infinitely many $m.$ Then, there is a subsequence $(y_{m_{k}})$
of $(y_{m})$ such that $y\succcurlyeq ^{\bullet }x\succcurlyeq y_{m_{k}}$,
so, by radiality, $d(x,y)\leq d(y,y_{m_{k}}),$ for all $k\geq 1$. Letting $%
k\rightarrow \infty $ here yields $x=y.$ Second, suppose $x\succcurlyeq
y_{m} $ holds for at most finitely many $m.$ Then, for some $M\in \mathbb{N}$%
, we have $x_{m}\succcurlyeq y_{m}\succcurlyeq ^{\bullet }x$ for all $m\geq
M $. Thus, by radiality, $d(x,y_{m})\leq d(x,x_{m})$ for all $m\geq M,$ so
letting $m\rightarrow \infty $ yields $x=y$ once again. We conclude that $%
y\succcurlyeq ^{\bullet }x$ implies $x=y,$ contradicting the reflexivity of $%
\succcurlyeq $. $\blacksquare $

\subsection{Partially Ordered Normed Spaces}

Let $X=(X,\left\Vert \cdot \right\Vert )$ be a normed (real) linear space.
We say that $(X,\succcurlyeq )$ is a \textit{partially ordered normed space}
if $\succcurlyeq $ is a \textit{vector order} on $X,$ that is, it is a
partial order on $X$ such that $x\succcurlyeq y$ implies $\lambda
x\succcurlyeq \lambda y$ for all $\lambda \geq 0$ and $x+z\succcurlyeq y+z$
for all $z\in X.$ If, in addition, $\succcurlyeq $ is a closed subset of $%
X\times X,$ we refer to $(X,\succcurlyeq )$ as a \textit{normed poset. }The
terms \textit{partially ordered Banach/Hilbert space} and \textit{%
Banach/Hilbert poset} are understood accordingly.

The \textit{positive cone }of any partially ordered normed space $%
(X,\succcurlyeq )$ is defined as $C(\succcurlyeq ):=\{x\in X:x\succcurlyeq 
\mathbf{0}\}$ where $\mathbf{0}$ denotes the origin of $X.$ As it contains $%
\mathbf{0},$ this set is always nonempty. As is well known, $C(\succcurlyeq
) $ is a pointed convex cone in $X$ such that $x\succcurlyeq y$ iff $x-y\in
C(\succcurlyeq )$.\footnote{%
As a reminder, we note that a \textit{convex cone} in a linear space $X$ is
a subset $C$ of $X$ such $\lambda C\subseteq C$ for all $\lambda \geq 0,$
and $C+C\subseteq C.$ Such a cone is said to be \textit{pointed} if $C\cap
-C=\{\mathbf{0}\}$.} (Conversely, for every pointed convex cone $C$ in $X,$
the binary relation $\succcurlyeq _{C}$ on $X$ with $x\succcurlyeq _{C}y$
iff $x-y\in C,$ is a vector order on $X$.) Obviously, $\succcurlyeq $ is
trivial iff $C(\succcurlyeq )=\{\mathbf{0}\}$. Moreover, if $(X,\succcurlyeq
)$ is a normed poset, $C(\succcurlyeq )$ is closed in $X$.

\subsection{Hadamard Posets}

Let $(X,\succcurlyeq )$ be a metric poset. If $X$ is uniquely geodesic, and $%
\succcurlyeq $ is \textit{geodesically hereditary}, that is, 
\begin{equation*}
x\succcurlyeq y\hspace{0.2in}\text{implies\hspace{0.2in}}x\succcurlyeq
z\succcurlyeq y\text{ for every }z\in \lbrack x,y]\text{,}
\end{equation*}%
where $[x,y]$ is the geodesic segment between $x$ and $y,$ we say that $%
(X,\succcurlyeq )$ is a \textit{geodesic poset}. If, in addition, $X$ is a
Hadamard space, that is, it is a complete CAT(0) space -- this entails
unique geodicity -- we say that $(X,\succcurlyeq )$ is a \textit{Hadamard
poset}. It is plain that every Hilbert poset is a Hadamard poset.

\subsection{Lipschitz Functions}

Let $X=(X,\succcurlyeq _{X})$ and $Y=(Y,\succcurlyeq _{X})$ be two partially
ordered metric spaces, and $f:X\rightarrow Y$ any function. For any real
number $K>0,$ we say that $f$ is $K$\textit{-Lipschitz} if $%
d_{Y}(f(x),f(y))\leq Kd_{X}(x,y)$ for every $x,y\in X$. (Here, and in what
follows, $d_{X}$ stands for the metric of $X,$ and similarly for $d_{Y}$.)
We say that $f$ is \textit{Lipschitz }if it is $K$-Lipschitz for some $K>0.$
For excellent treatments of the general theory of Lipschitz functions, see 
\cite{BeerBook}, \cite{C-M-N}, and \cite{Weaver}.

We denote the set of all $K$-Lipschitz maps from $X$ to $Y$ as Lip$%
_{K}(X,Y), $ but write Lip$_{K}(X)$ for Lip$_{K}(X,\mathbb{R}).$ In turn,
the sets of all order-preserving members of Lip$_{K}(X,Y)$ and Lip$_{K}(X)$
are denoted as Lip$_{K,\uparrow }(X,Y)$ and Lip$_{K,\uparrow }(X),$
respectively.

\subsection{The Monotone Lipschitz Extension Property}

Let $X=(X,\succcurlyeq _{X})$ and $Y=(Y,\succcurlyeq _{Y})$ be two partially
ordered metric spaces. We say that $(X,Y)$ has the \textit{monotone
Lipschitz extension property }if for every nonempty $S\subseteq X$ and $f\in 
$ Lip$_{1,\uparrow }(S,Y),$ there exists an $F\in $ Lip$_{1,\uparrow }(X,Y)$
with $F|_{S}=f.$ In this terminology, the classical \textit{McShane-Whitney
extension theorem} says: Provided that $\succcurlyeq _{X}$ is trivial, $(X,%
\mathbb{R})$ has the monotone Lipschitz extension property. Ok \cite{Ok} has
generalized this theorem as: $(X,\mathbb{R})$ has the monotone Lipschitz
extension property if, and only if, $\succcurlyeq _{X}$ is radial.

Again using this jargon, we can state \textit{Kirszbraun's theorem} as
follows: Provided that $X$ and $Y$ are Hilbert spaces, and $\succcurlyeq
_{X} $ is trivial, $(X,Y)$ has the monotone Lipschitz extension property.
This is a famous generalization of the McShane-Whitney theorem. The primary
objective of the present paper is to figure the extent to which one can
relax the triviality requirement on $\succcurlyeq _{X}$ in this statement.

By the \textit{monotone Lipschitz extension problem} relative to two
partially ordered metric spaces $X$ and $Y$, we understand determining if $%
(X,Y)$ has the monotone Lipschitz extension property.

\section{A Possibility Theorem}

Let us first consider the monotone Lipschitz extension problem in a very
specific case, namely, when $X=\mathbb{R}$. In that setup, the partial order
on $X$ is either trivial or linear. In the former case all maps that are
defined on any subset of $X$ are automatically order-preserving (regardless
of the codomain), so the problem reduces to a standard Lipschitz extension
problem. (As the proof of Theorem 2 shows, such an extension is indeed
possible, at least when $Y$ is a Banach space.) In the latter case too all
is well.

{\small \medskip }

{\small \medskip }

\noindent \textsc{Theorem 2.} \textsl{Let }$Y=(Y,\succcurlyeq )$ \textsl{be
a Banach poset. Then, }$(\mathbb{R},Y)$ \textsl{has the monotone Lipschitz
extension property.}

{\small \medskip }

{\small \medskip }

We can prove this by affine interpolation. For brevity, we will only sketch
the basic argument. Let $S$ be any nonempty set of the reals, and take any
order-preserving 1-Lipschitz $f:S\rightarrow Y$. For any $t_{1},t_{2},...\in
S$ with $t_{m}\rightarrow t,$ the Lipschitz property readily entails that $%
(f(t_{m}))$ is Cauchy in $Y,$ so as $Y$ is Banach, $(f(t_{m}))$ converges in 
$Y$. (Again by the Lipschitz property, the limit depends only on $t,$ not
the sequence $(t_{m}).$) We may then define $F:$ cl$(S)\rightarrow Y$ by $%
F(t):=\lim f(t_{m}),$ where $(t_{m})$ is any sequence in $S$ with $%
t_{m}\rightarrow t$. It is plain that $F$ is 1-Lipschitz. The closedness of $%
\succcurlyeq $ ensures that $F$ is order-preserving as well.

It is thus without loss of generality to assume $S$ is closed. If $S=\mathbb{%
R},$ there is nothing to prove, so we take $S$ as a nonempty closed proper
subset of $\mathbb{R}$. Then, there exists a countable collection $\mathcal{I
}$ of disjoint nonempty open intervals with $\mathbb{R}\backslash S=\bigcup 
\mathcal{I}$. For any $(a,b)\in \mathcal{I}$ with $a,b\in \mathbb{R},$ we
define 
\begin{equation*}
F(t):=f(a)+\frac{t-a}{b-a}(f(b)-f(a)),\,\,\,\,\,a<t<b.
\end{equation*}%
If $S$ has a minimum element, say, $a_{\ast }$, then $(-\infty ,a_{\ast
})\in \mathcal{I},$ and we set $F(t):=f(a_{\ast })$ for all $-\infty
<t<a_{\ast }$. Dually, if $S$ has a maximum element, say, $a^{\ast }$, then $%
(a^{\ast },\infty )\in \mathcal{I},$ and we set $F(t):=f(a^{\ast })$ for all 
$a^{\ast }<t<\infty $. This defines the $Y$-valued map $F$ everywhere on $%
\mathbb{R}$. It is routine to verify that $F$ is order-preserving and
1-Lipschitz.

\section{Necessity of Radiality}

It turns out that, under fairly general conditions, a pair of a partially
ordered metric spaces $X$ and $Y$ cannot have the monotone Lipschitz
extension property\textit{\ }unless $X$ is radially ordered. We will use
this fact as a springboard for the main results of this paper.

{\small \medskip }

{\small \medskip }

\noindent \textsc{Theorem 3.} \textsl{Let }$X=(X,\succcurlyeq _{X})$\textsl{%
\ and }$Y=(Y,\succcurlyeq _{Y})$ \textsl{be two partially ordered metric
spaces. Assume that there exists a geodesic ray }$\sigma :[0,\infty
)\rightarrow Y$ \textsl{such that}

(a) $\sigma (s)\succcurlyeq _{Y}\sigma (t)$ \textsl{whenever }$s\geq t\geq
0; $ \textsl{and}

(b) \textsl{the Busemann function }$B_{\sigma }$\textsl{\ associated with }$%
\sigma $\textsl{\ is} $\succcurlyeq _{Y}$\textsl{-decreasing.\footnote{%
We recall that $B_{\sigma }:Y\rightarrow \mathbb{R}$ is defined as $%
B_{\sigma }(a):=\lim_{t\rightarrow \infty }(d_{Y}(a,\sigma (t))-t).$ It is
well-known that this function is well-defined, 1-Lipschitz, and satisfies $%
B_{\sigma }(\sigma (t))=-t$ for all $t\geq 0$. As $Y$ is CAT(0), it is
convex as well. (See, for instance, Chapter 12 of \cite{Papa}.)}}

\noindent \textsl{If }$(X,Y)$\textsl{\ has the monotone Lipschitz extension
property, then }$\succcurlyeq _{X}$\textsl{\ is radial.}

{\small \medskip }

\noindent \textit{Proof.} Suppose $(X,Y)$\textsl{\ }has the monotone
Lipschitz extension property, but $\succcurlyeq _{X}$\textsl{\ }is not
radial. Then, there exist $x,y,z\in X$ such that either%
\begin{equation}
x\succcurlyeq _{X}^{\bullet }y\succ _{X}z\,\,\,\,\text{and}%
\,\,\,\,d_{X}(x,z)<d_{X}(x,y)  \label{birr}
\end{equation}%
or%
\begin{equation}
x\succ _{X}y\succcurlyeq _{X}^{\bullet }z\,\,\,\,\text{and}%
\,\,\,\,d_{X}(x,z)<d_{X}(y,z)\text{.}  \label{ikii}
\end{equation}%
We only consider the case (\ref{birr}), as the argument for (\ref{ikii}) is
analogous. Define $f:\{x,y\}\rightarrow Y$ by $f(x):=\sigma (d_{X}(x,y))$
and $f(y):=\sigma (0).$ As $d_{Y}(f(x),f(y))=\left\vert
d_{X}(x,y)-0\right\vert =d_{X}(x,y),$ it is plain that $f$ is 1-Lipschitz.
In addition, by hypothesis (a), $f(x)\succcurlyeq _{Y}f(y),$ so $f$ is
order-preserving. Then, by the monotone Lipschitz extension property, there
is an $F:X\rightarrow Y$ with $F(x):=\sigma (d_{X}(x,y))$ and $F(y):=\sigma
(0).$ Since $y\succ _{X}z,$ we have $F(y)\succcurlyeq _{Y}F(z).$ Then, by
hypothesis (b), 
\begin{equation}
0=-B_{\sigma }(\sigma (0))=-B_{\sigma }(F(y))\geq -B_{\sigma }(F(z))\text{.}
\label{uc}
\end{equation}%
Let us define $\varphi :=-B_{\sigma }\circ F$ so that $\varphi
(x)=-B_{\sigma }(\sigma (d_{X}(x,y)))=d_{X}(x,y)>0,$ while by (\ref{uc}), $%
0=\varphi (y)\geq \varphi (z).$ In particular, $\left\vert \varphi
(x)-\varphi (z)\right\vert =\varphi (x)-\varphi (z).$ On the other hand, as
both $-B_{\sigma }$ and $F$ are 1-Lipschitz, so is $\varphi :=-B_{\sigma
}\circ F$. Therefore, by (\ref{birr}), 
\begin{equation*}
d_{X}(x,y)-\varphi (z)=\varphi (x)-\varphi (z)\leq d_{X}(x,z)<d_{X}(x,y),
\end{equation*}%
that is, $\varphi (z)>0,$ a contradiction. $\blacksquare $

\section{Radial Orders on Connected Spaces}

Theorem 3 witnesses the importance of radiality for the monotone Lipschitz
extension problem. We are thus naturally led to investigate how demanding
the radiality property is. While \cite{Ok} provides quite a few examples of
radially ordered metric spaces, those are not particularly suitable for the
present setup, as they are all disconnected. In this section we will show
that this is not accidental. It turns out that there is a considerable
tension between connectedness of a metric space and its radial orderability.

{\small \medskip }

{\small \medskip }

\noindent \textsc{Theorem 4.} \textsl{Let }$X=(X,\succcurlyeq )$\textsl{\ be
a radially ordered metric space. If }$X$ \textsl{is connected, then }$%
\succcurlyeq $ \textsl{is either trivial or total.}

{\small \medskip }

\noindent \textit{Proof.} Assume that $X$ is connected, and take an
arbitrary $x\in X$. We define $x^{\uparrow }:=\{y\in X:y\succcurlyeq x\}$
and $x^{\downarrow }:=\{y\in X:x\succcurlyeq y\}.$ As $\succcurlyeq $ is
closed (Lemma 1), it is plain that these are closed subsets of $X$. We claim
that both $x^{\uparrow }\backslash \{x\}$ and $x^{\downarrow }\backslash
\{x\}$ are open in $X.$ To see this, take any $y\in x^{\uparrow },$ distinct
from $x$, and towards a contradiction, suppose $y$ is in the boundary of $%
x^{\uparrow }\backslash \{x\}$. Then, there is a sequence $(y_{m})$ in $%
X\backslash x^{\uparrow }$ with $y_{m}\rightarrow y.$ Consequently, $y\succ
x\succcurlyeq ^{\bullet }y_{m},$ so by radiality, $d(x,y_{m})\leq
d(y,y_{m}), $ for all $m.$ Then, letting $m\rightarrow \infty $ yields $x=y,$
contradicting $y$ being distinct from $x.$ This proves that $x^{\uparrow
}\backslash \{x\}$ is open in $X.$ Openness of $x^{\downarrow }\backslash
\{x\}$ is proved analogously. As $x$ was arbitrarily chosen in $X,$ and $X$
is connected, Schmeidler's theorem (of \cite{Schmeidler}) yields our claim. $%
\blacksquare $

{\small \medskip }

{\small \medskip }

The second alternative envisaged in Theorem 4 fails in myriad cases. To wit:

{\small \medskip }

{\small \medskip }

\noindent \textsc{Theorem 5.} \textsl{Let }$X$\textsl{\ be a connected
metric space that contains a copy of the circle }$\mathbb{S}^{1}$\textsl{.
Then, the only radial order on }$X$ \textsl{is the trivial order.}

{\small \medskip }

\noindent \textit{Proof.} Let $\succcurlyeq $ be a radial order on $X$. By
Theorem 4, $\succcurlyeq $ is either total or trivial. To derive a
contradiction, suppose the former case holds. Then, where $f:\mathbb{S}%
^{1}\rightarrow X$ is an embedding, we define the binary relation $%
\trianglerighteq $ on $\mathbb{S}^{1}$ by $a\trianglerighteq b$ iff $%
f(a)\succcurlyeq f(b).$ As $\succcurlyeq $ is closed (Lemma 1), $%
\trianglerighteq $ is a closed linear order $\trianglerighteq $ on $\mathbb{S%
}^{1}$. Without loss of generality, we may assume $e^{\pi i}\trianglerighteq
1.$ By closedness of $\trianglerighteq $, the set $\{\theta \in \lbrack
0,\pi ]:e^{(\theta +\pi )i}\trianglerighteq e^{\theta i}\}$ is a closed,
hence compact, subset of $\mathbb{R}$. It thus contains its supremum, say, $%
\theta ^{\ast }$. Clearly, $\theta ^{\ast }<\pi $ because otherwise, $%
1=e^{2\pi i}\trianglerighteq e^{\pi i},$ contradicting the antisymmetry of $%
\trianglerighteq $. It follows that 
\begin{equation*}
e^{(\theta ^{\ast }+\pi )i}\trianglerighteq e^{\theta ^{\ast }i}\text{
\thinspace \thinspace and \thinspace \thinspace }e^{(\theta ^{\ast
}+\varepsilon )i}\trianglerighteq e^{(\theta ^{\ast }+\varepsilon +\pi )i}
\end{equation*}%
for arbitrarily small $\varepsilon >0.$ This contradicts $\trianglerighteq $
being closed. $\blacksquare $

{\small \medskip }

{\small \medskip }

\noindent \textsc{Remark.} This proof is based on the observation that there
is no closed total order on any connected metric space that contains a
closed simple loop. The latter condition in this observation, and hence in
Theorem 5, can be relaxed quite a bit. See, for instance, \cite{N-O} and 
\cite{vanMill}.

{\small \medskip }

{\small \medskip }

We now extract some concrete cases out of Theorem 5. First, let us agree to
call a geodesic space $X$ \textit{nonbranching} if for any geodesics $\sigma
:[0,1]\rightarrow X$ and $\sigma ^{\prime }:[0,1]\rightarrow X$ with $\sigma
|_{[0,\varepsilon )}=\sigma ^{\prime }|_{[0,\varepsilon )}$ for some $%
\varepsilon >0,$ we have $\sigma =\sigma ^{\prime }.$ (For instance, every
strictly convex normed space is nonbranching, but an $\mathbb{R}$-tree is
not.) As usual, if $X$ is uniquely geodesic, we denote the geodesic segment
between any two points $x$ and $y$ in $X$ by $[x,y],$ and set $%
(x,y):=[x,y]\backslash \{x,y\}.$

{\small \medskip }

{\small \medskip }

\noindent \textsc{Corollary 6.} \textsl{Let }$X$\textsl{\ be a nonbranching
and uniquely geodesic space. Then, either }$X$ \textsl{is isometric to an
interval, or the only radial order on }$X$ \textsl{is the trivial order.}

{\small \medskip }

\noindent \textit{Proof.} We distinguish between two cases. First, assume
that every distinct three points in $X$ are collinear (in the sense that one
of these points lies in the geodesic segment between the other two points).
Fix any distinct $a,b\in X,$ and define $\tau :X\rightarrow \mathbb{R}$ by 
\begin{equation*}
\tau (x):=\left\{ 
\begin{array}{ll}
d(a,x), & \text{if }x\in \lbrack a,b]\text{ or }b\in \lbrack a,x], \\ 
-d(a,x), & \text{if }a\in \lbrack x,b]\text{.}%
\end{array}%
\right.
\end{equation*}%
Now take any distinct $x,y\in X.$ Without loss of generality, suppose $\tau
(y)\geq \tau (x).$ If $\tau (x)\geq 0,$ then $x\in \lbrack a,y],$ so $%
d(a,x)+d(x,y)=d(a,y),$ whence $\tau (y)-\tau (x)=d(a,y)-d(a,x)=d(x,y).$
Moreover, if both $\tau (y)$ and $\tau (x)$ are negative numbers, the same
argument applies. Finally, suppose $\tau (x)<0<\tau (y)$. In this case, $%
a\in \lbrack x,y],$ so $d(x,a)+d(a,y)=d(x,y),$ so $\tau (y)-\tau
(x)=d(a,y)-(-d(a,x))=d(x,y).$ We conclude that $d(x,y)=\left\vert \tau
(x)-\tau (y)\right\vert $ for all $x,y\in X.$ It follows that $\tau $ is an
isometry from $X$ into $\mathbb{R}$. As $X$ is geodesic, hence connected, $%
\tau (X)$ is an interval. Thus, $X$ is isometric to an interval.

Finally, suppose there exist distinct points $a,b,c\in X$ that are not
collinear. In this case $(a,b),$ $(b,c)$ and $(a,c)$ must be pairwise
disjoint. (For instance, if $d$ belongs to $(a,b)\cap (a,c),$ we must have $%
[a,d]\subseteq \lbrack a,b]$ and $[a,d]\subseteq \lbrack a,c]$ by unique
geodicity, but as $b\neq c,$ this contradicts $X$ being nonbranching.) It
follows that $S:=[a,b]\cup \lbrack b,c]\cup \lbrack a,c]$ is a closed simple
loop in $X$ which is, of course, homeomorphic to $\mathbb{S}^{1}$. By
Theorem 5, therefore, the only radial order on $X$ is the trivial one. $%
\blacksquare $

{\small \medskip }

{\small \medskip }

The following is a special case of Corollary 6.

{\small \medskip }

{\small \medskip }

\noindent \textsc{Corollary 7.} \textsl{Let }$X=(X,\succcurlyeq )$\textsl{\
be a radially ordered normed space. Then, either }$\dim X=1$\textsl{\ or }$%
\succcurlyeq $\textsl{\ is trivial.}

{\small \medskip }

{\small \medskip }

\noindent \textsc{Remark.} The no-extension theorem we have stated in the
Introduction is based on Corollary 7. It may thus be useful to give a
direct, alternative proof for this result that does not rely on Schmeidler's
theorem. Suppose $\dim X\geq 2,$ and to derive a contradiction, assume $%
\succcurlyeq $ is not trivial. Put $C:=C(\succcurlyeq )$, the positive cone
of $X$. We begin with showing that $C$ has a nonzero boundary point. Suppose
otherwise, so $C\backslash \{\mathbf{0}\}$ is an open set. As $\succcurlyeq $
is radial, it is closed (Lemma 1), and this entails $C$ is closed. It
follows that $C\backslash \{\mathbf{0}\}$ is a clopen set in $X\backslash \{%
\mathbf{0}\}.$ But, as $\dim X\geq 2,$ $X\backslash \{\mathbf{0}\}$ is
path-connected, hence connected. Thus, either $C\backslash \{\mathbf{0}%
\}=\varnothing $ or $C\backslash \{\mathbf{0}\}=X\backslash \{\mathbf{0}\}$.
In the former case we get $C=\{\mathbf{0\}},$ contradicting the
nontriviality of $\succcurlyeq ,$ and in the latter case $C=X,$
contradicting the antisymmetry of $\succcurlyeq $. We conclude that there is
a nonzero point $y$ in the boundary of $C$. Let us now pick any sequence $%
(y_{m})$ in $X\backslash C$ with $y_{m}\rightarrow y$. Then, $-y_{m}\in
X\backslash -C$ and $-y\in -C,$ that is, $-y_{m}\succcurlyeq ^{\bullet }%
\mathbf{0}\succcurlyeq -y$, for all $m\geq 1$. By radiality, then, $%
\left\Vert -y_{m}\right\Vert \leq \left\Vert -y_{m}+y\right\Vert $ for all $%
m\geq 1.$ Letting $m\rightarrow \infty $ yields $\left\Vert y\right\Vert
\leq 0,$ that is, $y=\mathbf{0},$ a contradiction. $\blacksquare $

\section{Impossibility Theorems}

Combining Theorem 3 and Corollary 6 yields the main result of this paper:

{\small \medskip }

{\small \medskip }

\noindent \textsc{Theorem 8.} \textsl{Let }$X=(X,\succcurlyeq _{X})$\textsl{%
\ be a partially ordered, nonbranching and uniquely geodesic space, and }$%
Y=(Y,\succcurlyeq _{Y})$\textsl{\ a partially ordered metric space. Suppose }%
$X$ \textsl{is not isometric to an interval, and there exists an
order-preserving geodesic ray }$\sigma :[0,\infty )\rightarrow Y$ \textsl{%
whose Busemann function is} $\succcurlyeq _{Y}$\textsl{-decreasing. If }$%
(X,Y)$\textsl{\ has the monotone Lipschitz extension property, then }$%
\succcurlyeq _{X}$\textsl{\ is trivial.}

{\small \medskip }

{\small \medskip }

The converse is false at this level of generality.\footnote{%
We thank Urs Lang for pointing this fact to us.} A pair of nonbranching
Hadamard spaces may well fail to have the Lipschitz extension property. For
a simple example, take the Poincare disk, and any equilateral triangle on
it. Map the corners of this triangle to an equilateral triangle (with same
side lengths) on the plane isometrically. This map cannot be extended to a
1-Lipschitz map on the Poincare disk.

{\small \medskip }

{\small \medskip }

\noindent \textsc{Corollary 9.} \textsl{Let }$X=(X,\succcurlyeq _{X})$%
\textsl{\ be a partially ordered Hilbert space with }$\dim X\geq 2$\textsl{,
and }$Y=(Y,\succcurlyeq _{Y})$\textsl{\ a Hadamard poset with an
order-preserving geodesic ray whose Busemann function is} $\succcurlyeq _{Y}$%
\textsl{-decreasing. Then, }$(X,Y)$\textsl{\ has the monotone Lipschitz
extension property if, and only if, }$\succcurlyeq _{X}$\textsl{\ is trivial.%
}

{\small \medskip }

\noindent \textit{Proof.} We only need to prove the sufficiency, but this is
a special case of Theorem A of Lang and Schroeder \cite{L-S}. $\blacksquare $

{\small \medskip }

{\small \medskip }

We next show that the conditions of Corollary 9 are satisfied by any Hilbert
poset $Y$.

{\small \medskip }

{\small \medskip }

\noindent \textsc{Proposition 10.} \textsl{Let }$Y=(Y,\succcurlyeq _{Y})$ 
\textsl{be a Hilbert poset. Then, there exists an order-preserving geodesic
ray }$\sigma :[0,\infty )\rightarrow Y$ \textsl{whose Busemann function is} $%
\succcurlyeq _{Y}$\textsl{-decreasing.}

{\small \medskip }

{\small \medskip }

The key to this result comes from convex analysis: For every Hilbert space $%
Y $ and a pointed closed and convex cone $C$ in $Y$, there is a positive
direction (relative to the partial order induced by $C$) that lies within
the dual of $C.$

{\small \medskip }

{\small \medskip }

\noindent \textsc{Lemma 11.} \textsl{Let }$Y=(Y,\left\langle \cdot ,\cdot
\right\rangle )$ \textsl{be a Hilbert space, and }$C$ \textsl{a pointed
closed and convex cone in }$Y$\textsl{\ with }$C\neq \{\mathbf{0}\}$\textsl{%
. Then, }$C\cap C^{\ast }\neq \{\mathbf{0}\}$\textsl{.}\footnote{%
We use the standard notation of convex analysis. In particular, $C^{\ast }$
stands for the \textit{dual cone} of $C,$ that is, $C^{\ast }=\{x\in
X:\left\langle x,y\right\rangle \geq 0$ for all $y\in C\}.$ The \textit{%
polar cone }of $C$ is $C^{\circ }:=-C^{\ast }$.}

{\small \medskip }

\noindent \textit{Proof.} Clearly, $C^{\ast }\neq \{\mathbf{0}\}$, because
otherwise $C^{\circ }=\{\mathbf{0}\},$ and by the bipolar theorem, $%
C=C^{\circ \circ }=Y,$ contradicting $C$ being pointed. Moreover, $C^{\ast
}\neq -C^{\ast },$ because otherwise $C^{\circ }$ is a linear subspace of $%
Y, $ and since $C=C^{\circ \circ },$ this implies $C$ is a linear subspace
of $Y,$ contradicting $C$ being pointed. It follows that there exists an $%
a\in C^{\ast }\backslash -C^{\ast }$.

Let $P_{C}$ denote the metric projection operator from $Y$ onto $C.$ By
Moreau's decomposition theorem, 
\begin{equation*}
a=P_{C}(a)+P_{C^{\circ }}(a).
\end{equation*}%
We set $b:=P_{C}(a).$ This vector is nonzero because otherwise, $%
a=P_{C^{\circ }}(a),$ so $a\in C^{\circ }=-C^{\ast },$ a contradiction.
Moreover, $a-b=P_{C^{\circ }}(a)\in C^{\circ },$ so $\left\langle
a-b,c\right\rangle \leq 0$ for all $c\in C.$ But, as $a\in C^{\ast },$ we
have $\left\langle a,c\right\rangle \geq 0$ for all $c\in C.$ It follows that%
\begin{equation*}
\left\langle b,c\right\rangle =\left\langle a,c\right\rangle -\left\langle
a-b,c\right\rangle \geq 0\,\,\,\,\text{for all }c\in C,
\end{equation*}%
that is, $b\in C^{\ast }.$ Conclusion, $b$ is a nonzero element of $C\cap
C^{\ast }$. $\blacksquare $

{\small \medskip }

{\small \medskip }

\textit{Proof of Proposition 10. }It is plain that $Y=(Y,\succcurlyeq _{Y})$
is an Hadamard poset. In addition, by Lemma 11, $C(\succcurlyeq _{Y})\cap
C(\succcurlyeq _{Y})^{\ast }\neq \{\mathbf{0}\},$ that is, there exists some 
$e\in Y$ such that $e\succ _{Y}\mathbf{0}$ and $\left\langle
e,a\right\rangle \geq 0$ for every $a\in Y$ with $a\succcurlyeq _{Y}\mathbf{0%
}$. Without loss of generality, we may take $\left\Vert e\right\Vert _{Y}=1.$
Consider the geodesic ray $\sigma :[0,\infty )\rightarrow Y$ with $\sigma
(t):=te.$ Since $e\succ _{Y}\mathbf{0},$ this ray is order preserving.
Moreover, for any $a\in Y$ and $t>0,$%
\begin{equation*}
\left\Vert a-\sigma (t)\right\Vert _{Y}-t=\frac{\left\Vert a-te\right\Vert
_{Y}^{2}-t^{2}}{\left\Vert a-te\right\Vert _{Y}+t}=\frac{\left\Vert
a\right\Vert _{Y}^{2}-2t\left\langle a,e\right\rangle }{\left\Vert
a-te\right\Vert _{Y}+t}=\frac{t^{-1}\left\Vert a\right\Vert
_{Y}^{2}-2\left\langle a,e\right\rangle }{\left\Vert t^{-1}a-e\right\Vert
_{Y}+1}.
\end{equation*}%
Letting $t\rightarrow \infty $ we find $B_{\sigma }(a)=-\left\langle
a,e\right\rangle $ for every $a\in Y,$ so $B_{\sigma }$ is indeed $%
\succcurlyeq _{Y}$-decreasing. $\blacksquare $

{\small \medskip }

{\small \medskip }

Combining Corollary 9 and Proposition 10 yields the no-extension theorem we
stated in the Introduction.

We conclude our exposition by noting that the generality afforded by Theorem
8 may be useful in other instances as well. The following simple examples
aim to demonstrate this.

{\small \medskip }

{\small \medskip }

\noindent \textsc{Example 1.} ($\mathbb{R}$\textsl{-tree-valued maps}) Let $%
Y $ be an $\mathbb{R}$-tree, that is, $Y$ is a uniquely geodesic metric
space such that for every $a,b,c\in Y$ with $[a,b]\cap \lbrack b,c]=\{b\},$
we have $[a,b]\cup \lbrack b,c]=[a,c].$\footnote{%
See \cite{A-O} and \cite{Evans} for nice introductions to the theory of $%
\mathbb{R}$-trees.} In addition, suppose the visual boundary of $Y$ is
nonempty, that is, there is at least one geodesic ray $\sigma :[0,\infty
)\rightarrow Y$. We define the binary relation $\succcurlyeq $ on $Y$ by 
\begin{equation*}
a\succcurlyeq b\text{\thinspace \thinspace \thinspace \thinspace
iff\thinspace \thinspace \thinspace \thinspace there is a }t\geq 0\text{
with }a\in \lbrack b,\sigma (t)].
\end{equation*}%
As $Y$ is uniquely geodesic, this is well-defined. It is readily checked
that $\succcurlyeq $ is a partial order, that is, $Y=(Y,\succcurlyeq )$ is a
poset. Obviously, $\sigma $ is order-preserving relative to this poset, and $%
a\succcurlyeq b$ and $c\in \lbrack a,b]$ implies $a\succcurlyeq
c\succcurlyeq b$.

Let us now compute the Busemann function associated with $\sigma $. Put $%
R:=\sigma ([0,\infty ))$, and take any $a\in Y.$ For every $t\geq 0,$ unique
geodicity of $Y$ entails that $[a,\sigma (t)]\cap R$ is a closed and
connected subset of $R$. Moreover, as $Y$ is an $\mathbb{R}$-tree, any two
geodesics that intersect do so in a (possibly degenerate) geodesic segment,
and once the geodesic from $a$ meets $R,$ it merges with it at a unique
point never to split again. It follows that there is a unique (hitting time) 
$t_{a}\geq 0$ such that $[a,\sigma (t)]\cap R$ equals $[\sigma
(t_{a}),\sigma (t)]$ for every $t\geq t_{a}$. Now put $m_{a}:=\sigma
(t_{a}). $ Then, $m_{a}\in \lbrack a,\sigma (t)],$ whence $d(a,\sigma
(t))=d(a,m_{a})+d(m_{a},\sigma (t)),$ for every $t\geq t_{a}$. Since $%
d(m_{a},\sigma (t))=t-t_{a},$ therefore, $d(a,\sigma
(t))-t=d(a,m_{a})-t_{a}, $ for every $t\geq t_{a}$. Conclusion:%
\begin{equation}
B_{\sigma }(a)=d(a,m_{a})-t_{a}\text{.}  \label{bb}
\end{equation}%
In particular, for any $a,b\in Y,$ setting $t_{ab}:=\max \{t_{a},t_{b}\},$
we have 
\begin{equation}
B_{\sigma }(b)-B_{\sigma }(a)=d(b,\sigma (t))-d(a,\sigma (t))  \label{bbb}
\end{equation}%
for every $t\geq t_{ab}$.

We will next characterize the partial order $\succcurlyeq $ in terms of $%
B_{\sigma }$. Take any $a,b\in Y$. Suppose $a\succcurlyeq b.$ Then, there is
a $T\geq 0$ with $a\in \lbrack b,\sigma (T)].$ So, $a\in \lbrack b,\sigma
(t)]$, whence $d(b,\sigma (t))=d(a,b)+d(a,\sigma (t)),$ for all $t\geq T.$
Substituting this in (\ref{bbb}) for, say, $t=\max \{T,t_{ab}\},$ we get $%
B_{\sigma }(b)-B_{\sigma }(a)=d(a,b).$ Conversely, assume $d(a,b)=B_{\sigma
}(b)-B_{\sigma }(a)$. Then, again by (\ref{bbb}), $d(b,\sigma
(t_{ab}))=d(a,b)+d(a,\sigma (t_{ab}))$. As $Y$ is uniquely geodesic, this
implies $a\in \lbrack b,\sigma (t_{ab})],$ that is, $a\succcurlyeq b.$ We
proved:%
\begin{equation}
a\succcurlyeq b\text{\thinspace \thinspace \thinspace \thinspace
iff\thinspace \thinspace \thinspace \thinspace }B_{\sigma }(b)-B_{\sigma
}(a)=d(a,b).  \label{ye}
\end{equation}%
As $(a,b)\mapsto B_{\sigma }(b)-B_{\sigma }(a)-d(a,b)$ is a continuous
function on $Y\times Y,$ it follows that $\succcurlyeq $ is closed, that is, 
$(Y,\succcurlyeq )$ is an Hadamard poset. Moreover, it is obvious from (\ref%
{ye}) that $B_{\sigma }$ is $\succcurlyeq $-decreasing. We may now apply
Theorem 8 to conclude: For any $n\geq 2,$ $(\mathbb{R}^{n},Y)$ does not have
the monotone Lipschitz extension property. $\square $

{\small \medskip }

{\small \medskip }

\noindent \textsc{Example 2.} (CAT(0)\textsl{-valued maps}) Let $Y$ be any
Hadamard space with a visual boundary (such as $\mathbb{H}^{n},$ $\mathbb{R}%
^{n}$, $\mathbb{H}^{n}\times \mathbb{R}^{n}$, etc.). In such a space, there
exists a geodesic ray $\sigma :[0,\infty )\rightarrow Y$. Define the binary
relation $\succcurlyeq $ on $Y$ as $a\succcurlyeq b$ iff $\{a,b\}\subseteq
\sigma ([0,\infty ))$ and $\sigma ^{-1}(x)\geq \sigma ^{-1}(y)$. It is
readily checked that this is a closed and geodesically hereditary partial
order on $Y$ relative to which $\sigma $ is order-preserving. As $B_{\sigma
}(\sigma (t))=-t$ for all $t\geq 0,$ it is also plain that $B_{\sigma }$ is $%
\succcurlyeq $-decreasing. By Theorem 8, therefore: For any $n\geq 2,$ $(%
\mathbb{R}^{n},Y)$ does not have the monotone Lipschitz extension property. $%
\square $

{\small \medskip }

{\small \medskip }

\noindent \textsc{Example 3.} ($\mathbb{H}^{n}$\textsl{-valued maps}) Let us
adopt the hyperboloid model for the hyperbolic $n$-space $\mathbb{H}^{n}$.
(That is, $\mathbb{H}^{n}$ is the metric space whose ground set is $\{a\in 
\mathbb{R}^{n+1}:\left\langle a,a\right\rangle _{L}=-1$ and $a_{n+1}>0\},$
and whose metric $d$ assigns to each $(a,b)\in \mathbb{H}^{n}\times \mathbb{H%
}^{n}$ the unique nonnegative number $d(a,b)$ with $\cosh
d(a,b)=-\left\langle a,b\right\rangle _{L}.$ Here $\left\langle \cdot ,\cdot
\right\rangle _{L}$ is the Lorentz bilinear form defined as $\left\langle
a,b\right\rangle _{L}:=\sum^{n}a_{i}b_{i}-a_{n+1}b_{n+1}$.) We define the
closed partial order $\succcurlyeq $ on $\mathbb{H}^{n}$ by $a\succcurlyeq b$
iff $a_{i}=b_{i}$ for each $i=1,...,n,$ and $a_{n+1}\geq b_{n+1}$, and set $%
Y=(\mathbb{H}^{n},\succcurlyeq ).$ Now consider the geodesic ray $\sigma
:[0,\infty )\rightarrow \mathbb{H}^{n}$ with $\sigma (t):=(0,...,0,e^{t}).$
It is routine to compute that $B_{\sigma }(a)=-\log a_{n+1}$ for every $a\in 
\mathbb{H}^{n}$. Thus, $B_{\sigma }$ is $\succcurlyeq $-decreasing. As the
other requirements of Theorem 8 are also easily met, we conclude: For any $%
n\geq 2,$ $(\mathbb{R}^{n},Y)$ does not have the monotone Lipschitz
extension property. $\square $

\section{Concluding Comments}

Let us restate the no-extension theorem of the Introduction in the common
jargon of the recent literature on Lipschitz extensions. Let $%
X=(X,\succcurlyeq _{X})$ and $Y=(Y,\succcurlyeq _{Y})$ be metric posets, and
for any $k\in \mathbb{N}\cup \{\infty \}$, let $\mathcal{X}_{k}$ denote the
collection of all $S\subseteq X$ with $0<\left\vert S\right\vert \leq k$.
For any nonempty $S\subseteq X$, define $e_{\uparrow }(X,S,Y):=\inf \{K\in
(0,\infty ]:$ for every $f\in $ Lip$_{1,\uparrow }(S,Y)$ there is an $F\in $
Lip$_{K,\uparrow }(X,Y)$ with $F|_{S}=f\},$ and for any $k\in \mathbb{N}\cup
\{\infty \}$, define%
\begin{equation*}
e_{k,\uparrow }(X,Y):=\sup_{S\in \mathcal{X}_{k}}e_{\uparrow }(X,S,Y).
\end{equation*}%
In the case where $\succcurlyeq _{X}$ is trivial (so Lip$_{1,\uparrow
}(S,Y)= $ Lip$_{1}(S,Y)$ for all $S\in \mathcal{X}$), we write $e(X,Y)$ for $%
e_{\uparrow ,\infty }(X,Y),$ and $e_{k}(X,Y)$ for $e_{k,\uparrow }(X,Y)$ for
any $k\in \mathbb{N}$.

Quite a bit is known about the quantities $e(X,Y)$ and $e_{k}(X,Y)$. For one
thing, the McShane-Whitney theorem says $e(X,\mathbb{R})=1$ -- in fact, $%
e(X,\ell _{\infty })=1$ -- whereas Kirszbraun's theorem says $e(X,Y)=1$ when 
$X$ and $Y$ are Hilbert spaces. In a famous paper, Johnson and Lindenstrauss 
\cite{J-L} have shown that $e_{k}(X,Y)\leq 2\sqrt{\log k}$ when $Y$ is a
Hilbert space. In the case where $Y$ is a Banach space, Johnson,
Lindenstrauss and Schlechtman \cite{JLS} have proved that $e_{k}(X,Y)=O(\log
k),$ and Lee and Naor \cite{Lee-Naor} have improved that estimate to $%
e_{k}(X,Y)=O(\frac{\log k}{\log \log k}).$ Remarkably, these fascinating
results are all obtained by probabilistic methods.\footnote{%
There are powerful estimates that bound $e_{k}(X,Y)$ from below as well. See 
\cite{Basso} to learn about the state of the art in this area of research.}

Very little is known in the order-theoretic setup, however. Ok \cite{Ok}
shows that $e_{\uparrow }(X,\mathbb{R})=e_{\uparrow }(X,\ell _{\infty })=1$
when $\succcurlyeq _{X}$ is radial, but we have found in this paper that
radiality is excessively demanding in the case of connected metric spaces.
We have $e_{\uparrow }(\mathbb{R},Y)=1$ whenever $Y$ is a Banach poset, but
this is little comfort. For instance, we even have $e_{\uparrow }(\mathbb{R}%
^{2},\mathbb{R})>1.$ More generally, the no-extension theorem of the
Introduction shows that when $X$ and $Y$ are Hilbert posets with $\dim X\geq
2,$ we have $e_{\uparrow }(X,Y)>1$. In fact, our proof of this result shows
that $e_{2,\uparrow }(X,Y)>1$ in this case. It should thus be quite
interesting to seek upper bounds for the quantity $e_{k,\uparrow }(X,Y)$
when $X$ and $Y$ are arbitrary Hilbert posets (with $\dim X\geq 2$), but
this will very much depend on the structure of the partial orders of $X$ and 
$Y$. For instance, as $e_{\uparrow }(X,\mathbb{R})=1$, we have $e_{\uparrow
}(X,\mathbb{R}^{n})\leq \sqrt{n}$ for all $n\in \mathbb{N}$ (where $X$ is
any metric poset), but we do not presently know if this is the best upper
bound, nor have good estimates for $e_{\uparrow ,k}(X,\mathbb{R}^{n})$.%
{\small \medskip }

\end{document}